\documentclass[final,1p,times]{elsarticle}
\usepackage{amssymb,color,ulem,amsmath}
\newcommand{\sign}{{\rm \hskip0.5pt sign \hskip1pt}}

\usepackage[usenames,dvipsnames]{xcolor}

\def\qed{\hfill $\square$}

\journal{}
\usepackage{amsthm,verbatim}
\theoremstyle{definition}

\newtheorem{thm}{Theorem}[section]
 
\newtheorem{rem}{Remark}[section]

\newtheorem{cor}{Corollary}[section]

\newtheorem{prop}{Proposition}[section]
\newproof{pf}{Proof}
\newcommand{\eps}{\varepsilon}
\begin{document}
\begin{frontmatter}

\author[mymainaddress]{Oleg Makarenkov\corref{mycorrespondingauthor}}
\cortext[mycorrespondingauthor]{Corresponding author}
\ead{makarenkov@utdallas.edu}

\author[mymainaddress]{Anthony Phung}

\address[mymainaddress]{Department of Mathematical Sciences, University of Texas at Dallas, 75080, TX, Richardson, USA}

\title{Dwell time for local stability of switched systems with application to non-spiking neuron models}
\begin{abstract} For switched systems that switch between distinct globally stable equilibria, we offer closed-form formulas that lock  oscillations in the required neighborhood of the equilibria. Motivated by non-spiking neuron models, the main focus of the paper is on the case of planar switched affine systems, where we use properties of nested cylinders coming from quadratic Lyapunov functions. In particular, for the first time ever, we use the dwell-time concept in order to give an explicit condition for non-spiking of linear neuron models with  periodically switching current. An extension to the general nonlinear case is also given.
\end{abstract}
\end{frontmatter}
\section{Introduction}\label{sec:int}

\noindent 
Dwell time is the lower bound on the time between successive discontinuities (switchings) 
of the piecewise constant function $u(t)$, which ensures that the corresponding switched system
\begin{equation}\label{ss}
   \dot x=f_{u(t)}(x),\quad  x\in\mathbb{R}^n,
\end{equation}
exhibits a required type of stability,  under the assumption that each of the subsystems 
\begin{equation}\label{ssi}
   \dot x=f_{u}(x),\quad u\in\mathbb{R},\ 
x\in\mathbb{R}^n,
\end{equation}
 possess a unique globally asymptotically stable equilibrium $x_u.$ Let  $V_u$ be some Lyapunov function of subsystem (\ref{ssi}) corresponding to $x_u$ and let $N_u^k$ be the neighborhood of $u$ given by
 \begin{equation}\label{Ni}
   N_u^k=   \left\{x:V_u(x)\le k \right\}.
 \end{equation}
 Extending the pioneering result by Alpcan-Basar \cite{alpcan} (see also Liberzon \cite[\S3.2.1]{liberzon}), the recent paper \cite{dorothy} by Dorothy-Chung  gives a formula for the dwell time $\tau_d$ which ensures that any solution of (\ref{ss}) with the initial condition $x(t_0)\in N_{u(t_0)}^k$ satisfies
 \begin{equation}\label{sat1}  
     x(t_i)\in N_{u(t_i)}^k,\quad  i\in\mathbb{N},
 \end{equation}
 as long as the successive discontinuities $t_1,t_2,...$ of the control signal $u(t)$ verify
 \begin{equation}\label{sat2}
    t_{i+1}-t_i\ge\tau_d,\quad\ i\in \mathbb{N}.
 \end{equation}
However, the results of \cite{dorothy} are formulated in general abstract settings and certain work is required to apply those results to particular problems. In the present paper we follow the strategy of \cite{dorothy} when addressing planar switched affine systems, but carry out an independent proof that allows us to get closed-form formulas for the dwell-time $\tau_d$ (i.e. formulas in terms of just coefficients of the affine subsystems). 

\vskip0.2cm

\noindent Relevant results have been recently obtained in Xu et al \cite{xu} for quasi-linear switched systems (\ref{ss}), but the dwell-time formula \cite{xu} is not fully explicit, as it involves the constant of the rate of decay of the matrix exponent of the homogeneous part of subsystems (\ref{ssi}).

\vskip0.2cm

\noindent Our research is motivated by an  application to non-spiking of linear neurons with a periodically switching current.
The model of a planar linear neuron reads as (Izhikevich \cite[\S8.1.1]{izh}, Hasselmo-Shay \cite{hasselmo}) 
\begin{equation}\label{eq1}
\begin{array}{l}
   \dot v=-g_p v+g_h h+I_{in}(t),\\
   \dot h=-mv-o_h h,
\end{array}  
\end{equation}
coupled with the reset law
\begin{equation}\label{reset}
 v(t+0)=v_R,\  h(t+0))=h_R(u(t-0)),\ \mbox{if}\ v(t)=v_{th},
\end{equation}
where $v$ is the neural cell membrane potential, $h$ is the recovery current,
$g_p$ is the rate of passive decay of membrane potential, $g_h$ is the rate of current induced depolarization of the cell, $m>0$ makes $h$ increasing when $v$ gets negative, $o_h$ is the current decay,
$I$ is a constant current which can switch on and switch off.  Though some neurons spike and reset according to $v_R$ and $h_R$ (when reach the threshold $v_{th}$) to propagate message, some others are capable to transmit information without spiking and are not supposed to ever reach the firing threshold $v_{th}$ (see e.g. Vich-Guillamon \cite{non-spiking1}, Chen et al \cite{non-spiking2}).
The present paper uses the dwell time concept in order to obtain conditions for the model (\ref{eq1})-(\ref{reset}) to never reach the firing threshold $v_R$, i.e. to ensure just subthreshold oscillations. The readers interested in the difference between subthreshold and spiking dynamics are referred to Coombes et al \cite{coombes}.

\vskip0.2cm

\noindent As for nonlinear switched systems with arbitrary Lyapunov functions $V_u$ (whose level curves are not necessary ellipses), we don't see how the strategy of 
\cite{dorothy} (that we use in the case of affine subsystems) can provide computable formulas for the dwell time. That is why we offer a different approach based on approximating the level curves of $V_u$ by inner and outer balls. Formally speaking, such an approach provides a more conservative bound for the dwell time compared to \cite{dorothy}, but the advantage of our approach is that it requires minimal computations.  
\vskip0.2cm

\noindent We would like to clarify that switched systems that switch between different equilibria upon crossing a threshold of phase space often stabilize to a point  on the threshold, which can be considered an equilibrium of switched system in a certain generalized sense (related to the notions of Filippov's sliding vector field and switched equilibrium), see e.g. Polyakov-Fridman \cite{fridman} and  
Bolzern-Spinelli \cite{bolzen}. The dwell-time concept deals with switched systems that switch in time, which is simpler to implement in practice, but which cannot provide converge to a required point of phase space (the time-based switching rule cannot sense the phase coordinate of the solution).

\vskip0.2cm

\noindent The paper is organized as follows. In the next section, we consider 2-dimensional switched affine systems of differential equations and give an explicit description of the trapping region that the solution of (\ref{ss}) with the initial condition  $x(t_0)\in N_{u(t_0)}^\eps$ belongs to when the switching instances $t_1,t_2,...$ satisfy (\ref{sat2}). Specifically, on top of (\ref{sat1}) we establish (Theorem~\ref{proplinear}) that $
   x(t)\in N^{k_i}_{u(t_i)},$ $t\in[t_{i-1},t_i],\ i\in\mathbb{N},$ where  $k_i$ is given by a explicit formula (\ref{ki}).  The proof is carried out by deriving a closed-form formula for the Lyapunov functions of affine subsystems (\ref{ssi}) and by constructing the ellipses $N^{k_i}_{u(t_i)}$ to contain and just touch $N^{k}_{u(t_{i-1})}$. Since our main goal is the linear neuron model (\ref{eq1}) with switched input $I_{in}$ we focus in section~2 on affine subsystems (\ref{ssi}) with $u$-independent homogeneous part, but we explain in Remark~\ref{rem1} how the proof extends to the case of $u$-dependent homogeneous parts. The result of Section~2 is then used  in Section~3 in order to locate  (Proposition~\ref{propneuron}) the trapping region of the linear neuron model (\ref{eq1}) with switching current $I_{in}$ and to give conditions for non-spiking. Furthermore, simulations of Section~3 document that the proposed estimate for the trapping region of neuron model (\ref{eq1}) is sharp (as well as the proposed formulas for the dwell time). Section~4 offers a possible method to obtain closed-form dwell-time formulas in the case where the level curves of the Lyapunov functions of subsystems (\ref{ssi})  are not ellipses. 
Conclusion and Acknowledgment sections conclude the paper.

\section{Dwell-time and local trapping region for planar switched affine systems}

\noindent In this section, we spot a situation where the strategy of  Dorothy-Chung \cite{dorothy} leads to  closed-form dwell-time formulas. Specifically, we consider the case where the subsystems (\ref{ssi}) have the form
\begin{equation}\label{ss_affine}
    \dot x =\left(\begin{array}{cc} 
       a & b\\ c & d
    \end{array}
    \right) x+B_u,
\end{equation}
with
\begin{equation}\label{abcd}
a b c d<0.
\end{equation} 


\noindent Observe, the matrix $\left(\begin{array}{cc} 
       a & b\\ c & d
    \end{array}
    \right)$ is Hurwitz, if (see Zhang et al \cite{automatica})
\begin{equation}\label{hur}
    ad-bc>0\quad\mbox{and}\quad a+d<0.
\end{equation}
\begin{prop}\label{propsimplify} If (\ref{abcd})-(\ref{hur}) hold, then 
    \begin{equation}\label{equilibria}
       V_u(x)=\sign(ac)(-c(x_1-x_{u,1})^2+b(x_2-x_{u,2})^2),\quad\mbox{where}\quad x_u=-\left(\begin{array}{cc} 
       a & b\\ c & d
    \end{array}
    \right)^{-1}B_u,
    \end{equation}
    is a positive definite Lyapunov function of the globally asymptotically stable system (\ref{ss_affine}).
\end{prop}

\noindent {\bf Proof.} We look for $V_u$ in the form
$$
   V_u(x)=x^TPx,
$$
where the $2\times 2$-matrix $P$  is the solution of the Lyapunov equation (see Vidyasagar \cite[Sec 5.4, Theorem 42]{77}, Khalil \cite[Theorem 3.6]{78})
$$
  2\sign(ac)\left(\begin{array}{cc} -ac & 0 \\ 0 & bd\end{array}\right)=\left(\begin{array}{cc} 
       a & c\\ b & d
    \end{array}
    \right)P+P\left(\begin{array}{cc} 
       a & b\\ c & d
    \end{array}
    \right).
$$
The conclusion is achieved by noticing that the required solution $P$ is given by
$$
   P=\sign(ac)\left(\begin{array}{cc} 
       -c & 0\\ 0 & b
    \end{array}
    \right).
$$
\qed


\begin{thm}\label{proplinear} Assume that  (\ref{abcd})-(\ref{hur}) hold.
Let $x$ be any solution of switched system (\ref{ss}) with the initial condition $x(t_0)\in N_{u(t_0)}^k$. If the successive discontinuities $t_1,t_2,...$ of the control signal $u(t)$ verify
 \begin{equation}\label{sat2new}
    t_{i}-t_{i-1}\ge\frac{1}{2\min\{|a|,|d|\}}\ln\left(
    \frac{k_i}{k}\right),\quad\ i\in \mathbb{N},
 \end{equation}
 where 
 \begin{equation}\label{ki}
 k_i=\left(\sqrt{k}+\sqrt{|c|\left(x_{u(t_{i}),1}-x_{u(t_{i-1}),1}\right)^2+|b|\left(x_{u(t_{i}),2}-x_{u(t_{i-1}),2}\right)^2}\right)^2,\quad i\in\mathbb{N},
 \end{equation}
 with
 the equilibria $x_u$ given by (\ref{equilibria}), then
\begin{equation}\label{dyn1}
   x(t_i)\in N_{u(t_i)}^k,\quad i\in\mathbb{N},
\end{equation}
and
\begin{equation}\label{dyn2}
   x(t)\in N^{k_i}_{u(t_i)},\quad t\in[t_{i-1},t_i],\ i\in\mathbb{N}.
\end{equation}
\end{thm}

\begin{figure}[h]\center
\includegraphics[scale=0.55]{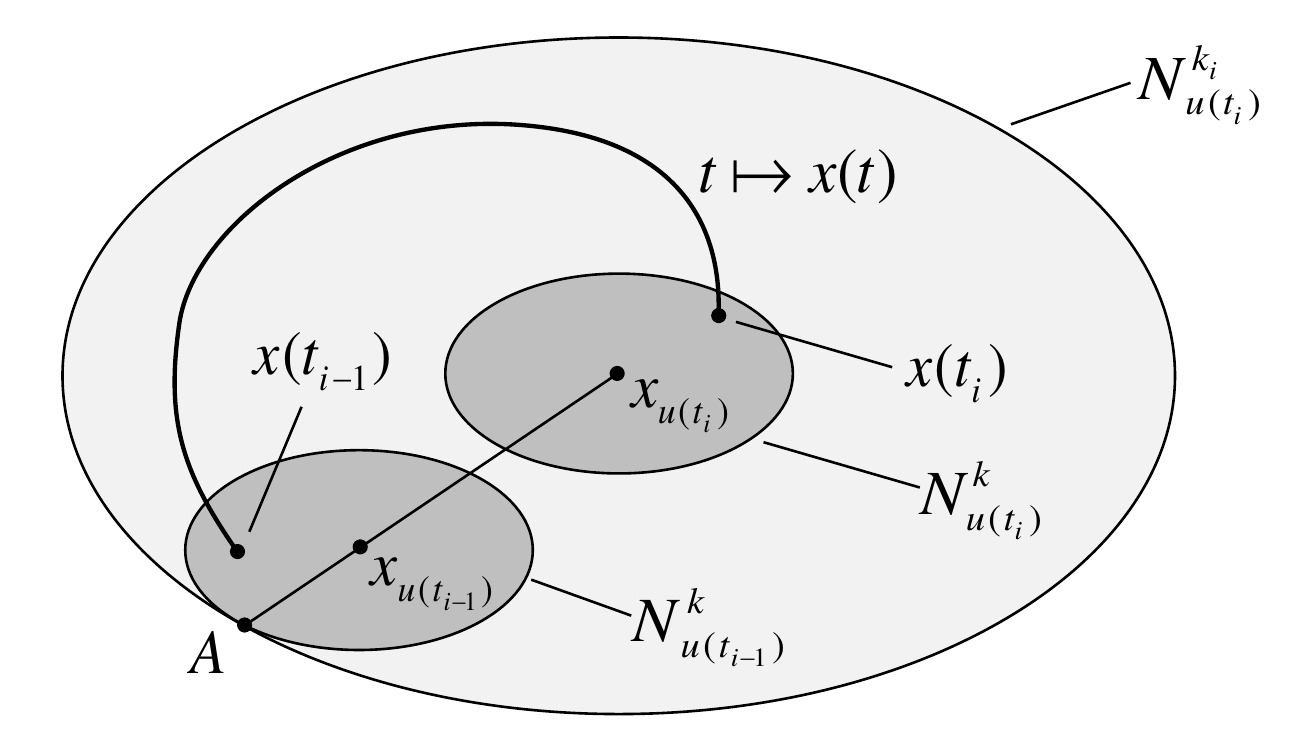}
\vskip-0.2cm
\caption{\footnotesize  The location of the ellipse $\partial N_{u(t_{i})}^{k_i}$ relative to ellipse $\partial N_{u(t_{i-1})}^{k}$ and the solution $x \mapsto x(t)$ of switched system (\ref{ss}) on the interval $[t_{i-1},t_i]$.} \label{figX}
\end{figure}

\vskip0.2cm 
\noindent The notations and conclusions of the theorem are illustrated in Fig.~\ref{figX}.

\noindent {\bf Proof.} 
Put
$$
   v_i(t)=V_{u(t_i)}(x(t)),\quad t\in(t_{i-1},t_i].
$$
Then
\begin{equation}\label{then1}
   \dot v_i(t)\le -\eps v_i(t),
\end{equation}
where $\eps>0$ is such a constant that
$$
\eps\le \frac{\sign(ac)\left(2a c (x_1-x_{u,1})^2-2bd(x_2-x_{u,2})^2\right)}{\sign(ac)\left(-c(x_1-x_{u,1})^2+b(x_2-x_{u,2})^2\right)}=2\frac{|a c| (x_1-x_{u,1})^2+|bd|(x_2-x_{u,2})^2}{|c|(x_1-x_{u,1})^2+|b|(x_2-x_{u,2})^2}.
$$
Letting $x_1-x_{u,1}=r\cos\phi$ and $x_2-x_{u,2}=r\sin\phi$, the right-hand-side of this inequality takes the form $$
    2\dfrac{|a c| \cos^2\phi+|bd|\sin^2\phi}{|c|\cos^2\phi+|b|\sin^2\phi}=:g(\phi).
$$ To find the best (i.e. maximal) possible value of $\eps$ we therefore compute the minimum of $g(\phi)$ on the interval $[0,\pi].$ We have
$$
  g'(\phi)=\dfrac{2|bc|(|a|-|d|)\sin\phi\cos\phi}{|c|\cos^2\phi+|b|\sin^2\phi}
$$
and so $g(\phi)$ has just one critical point $\phi_0=\pi/2$ on $(0,\pi).$ Therefore, 
\begin{equation}\label{then2}
  \eps=\min\limits_{\phi\in[0,\pi]}g(\phi)=\min\left\{g(0),g(\pi/2),g(\pi)\right\}=2\min\{|a|,|d|\}.
\end{equation}
Let us fix $i\in\mathbb{N}.$  Assuming that  $x(t_{i-1})\in N_{u(t_{i-1})}^k$ is established, we now use
(\ref{then1})-(\ref{then2}) in order to prove that
$x(t_{i})\in N_{u(t_{i})}^k$, i.e. to prove that 
$v_i(t_i)\le k.$
 Specifically, we are going to find  $k_i>0$ satisfying 
 \begin{equation}\label{NN}
 N_{u(t_{i-1})}^k\subset  N_{u(t_{i})}^{k_i}
 \end{equation}  and prove that 
\begin{equation}\label{provethat}
k_i{\rm e}^{-\eps_{u(t_i)}(t_i-t_{i-1})}\le k
\end{equation} to have
$$
   v_i(t_i)\le v_i(t_{i-1}){\rm e}^{-\eps_{u(t_i)}(t_i-t_{i-1})}\le  k_i{\rm e}^{-\eps_{u(t_i)}(t_i-t_{i-1})}\le k.
$$
Note, that  the boundary $\partial N_u^k$ of $N_u^k$ is given by
$$
   \partial N_u^k = \left\{x\in\mathbb{R}^2:|c|(x_1-x_{u,1})^2+|b|(x_2-x_{u,2})^2=k\right\}.
$$  To find $k_i>0$ satisfying (\ref{NN}) we construct the ellipse $\partial N_{u(t_{i})}^{k_i}$ to touch the ellipse $\partial N_{u(t_{i-1})}^{k}$, see Fig.~\ref{figX}. 
Let $A\in\mathbb{R}^2$ be the point where the two ellipses touch one another. 
Expressing the point $A$ in the polar coordinates of the ellipses $\partial N_{u(t_{i-1})}^{k}$ and $\partial N_{u(t_{i})}^{k_i}$ we get 
\begin{equation}\label{tmp1}
 \begin{array}{rcl}
   x_1-x_{u(t_{i-1}),1}&=&\sqrt{\frac{k}{{|c|}}}\cos\phi,\\
   x_2-x_{u(t_{i-1}),2}&=&\sqrt{\frac{k}{{|b|}}}\sin\phi,
 \end{array} \qquad\mbox{and}\qquad  \begin{array}{rcl}
   x_1-x_{u(t_{i}),1}&=&\sqrt{\frac{k_i}{{|c|}}}\cos\bar\phi,\\
   x_2-x_{u(t_{i}),2}&=&\sqrt{\frac{k_i}{{|b|}}}\sin\bar\phi.
 \end{array}
\end{equation}
The property of the derivative of 
 the curve $\partial N_{u(t_{i-1})}^{k}$ at $A$ to be parallel to the derivative of the curve $\partial N_{u(t_{i})}^{k_i}$ at $A$ leads to $\phi=\bar\phi.$ Excluding in (\ref{tmp1}) the unknowns $x_1$ and $x_2$ we get
 \begin{equation}\label{2root}
    x_{u(t_{i}),1}-x_{u(t_{i-1}),1}=\frac{\sqrt{k}-\sqrt{k_i}}{\sqrt{|c|}
}\cos\phi,\qquad x_{u(t_{i}),2}-x_{u(t_{i-1}),2}=\frac{\sqrt{k}-\sqrt{k_i}}{\sqrt{|b|}
}\sin\phi,
 \end{equation}
 which yields (\ref{ki}). Combining (\ref{ki}) with  (\ref{then2}), the inequality (\ref{provethat}) takes form of assumption (\ref{sat2new}) and so (\ref{provethat})  holds true. The proof of the theorem is complete.
\qed
 
 \vskip0.2cm
 
\begin{rem}\label{rem1} When system (\ref{ss_affine}) has the form
\begin{equation}\label{ssu}
    \dot x =\left(\begin{array}{cc} 
       a_u & b_u\\ c_u & d_u
    \end{array}
    \right) x+B_u,
\end{equation}
the 
property of the derivative of 
 the curve $\partial N_{u(t_{i-1})}^{k}$ at $A$ to be parallel to the derivative of the curve $\partial N_{u(t_{i})}^{k_i}$ at $A$ (see the proof of Theorem~\ref{proplinear}) no longer leads to $\phi=\bar\phi,$ but gives a relation
$$
   \frac{\sqrt{|c_{u(t_{i-1})}|}\cos\phi}{\sqrt{|b_{u(t_{i-1})}|}\sin\phi}=\frac{\sqrt{|c_{u(t_{i})}|}\cos\bar\phi}{\sqrt{|b_{u(t_{i})}|}\sin\bar\phi}.
$$
This relation needs to be used in order to eliminate the unknowns $x_{0,1},$ $x_{0,2}$, $\phi$, and $\bar\phi$ from the system
$$
 \begin{array}{rcl}
   x_1-x_{u(t_{i-1}),1}&=&\sqrt{\frac{k}{{|c_{u(t_{i-1})}|}}}\cos\phi,\\
   x_2-x_{u(t_{i-1}),2}&=&\sqrt{\frac{k}{{|b_{u(t_{i-1})}|}}}\sin\phi,
 \end{array} \qquad\mbox{and}\qquad  \begin{array}{rcl}
   x_1-x_{u(t_{i}),1}&=&\sqrt{\frac{ k_i}{{|c_{u(t_{i})}|}}}\cos\bar\phi,\\
   x_2-x_{u(t_{i}),2}&=&\sqrt{\frac{k_i}{{|b_{u(t_{i})}|}}}\sin\bar\phi,
 \end{array}
$$
which is the analogue of (\ref{tmp1}) when (\ref{ssu}) is considered in place of (\ref{ss_affine}). As a consequence, one gets a relation between $k$ and $k_i,$ which will replace (\ref{ki}) in the formulation of Theorem~\ref{proplinear} for planar switched affine systems of form (\ref{ssu}).
\end{rem}

\begin{rem} Further to Remark~\ref{rem1}, if subsystem (\ref{ssi}) with $u=u(t_i)$ is unstable, then the solution $x$ with the initial condition $x(t_{i-1})\in N^k_{u(t_{i-1})}$ never reaches $N^k_{u(t_i)}$ as long as $u(t)$ stays equal $u(t_i)$. In contrast, the trajectory $t\mapsto x(t)$ will go away from $N^k_{u(t_i)}$. In this case, one can use the ellipses of Lyapunov function $V_{u(t_i)}$ in order to evaluate how far will the trajectory $x$ deviate from $x_{u(t_i)}$ during the time $t_i-t_{i-1}$ (which will be then used  to allow extra time for the convergence to the stable equilibrium $x_{u(t_{i+1})}$ during $[t_{i},t_{i+1}]$). In this type of analysis one will need to construct such an ellipse $N_{u(t_i)}^{k_i}$, which just touches $N^k_{u(t_{i-1})}$, but don't cover it.  The respective $k_i=k_i^{unstab}$ will be the smaller root of (\ref{2root}), while the $k_i$ given by (\ref{ki}) was the largest root of (\ref{2root}). In this way, Theorem~\ref{proplinear} can be extended to the case  where some of the subsystems (\ref{ssi}) are unstable, complimenting the  results by Dorothy-Chung \cite{dorothy}, Zhai et al \cite{michel}, Li et al \cite{li}.
\end{rem}

\begin{rem} It is possible to extend Theorem~\ref{proplinear} to the multi-dimensional case under the assumption that the level curves of $V_u$ are still ellipsoids.
\end{rem}

\begin{rem} Note, Theorem~\ref{proplinear} is equally applicable in the case of a finite number of switchings $t_1,t_2, ...$ .
\end{rem}

\begin{rem} The dwell time requirement (\ref{sat2new}) can be weakened to 
$$    t_{i}-t_{i-1}\ge\frac{1}{2\min\{|a|,|d|\}}\ln\left(1+\frac{\max\left\{\sqrt{|c|},\sqrt{|b|}\right\}\cdot\|x_{u(t_i)}-x_{u(t_{i-1})}\|}{\sqrt{k}}
 \right),\quad\ i\in \mathbb{N}.
$$
\end{rem}

\section{Application to non-spiking neuron models} 

\noindent In this section we apply the earlier results to the planar linear system (\ref{eq1}) assuming that the current $I_{in}$ is changing according to the law
\begin{equation}\label{controlI}
   I_{in}(t)=\left\{\begin{array}{ll}
     I, & t\in(0,T_I],\\
     0, & t\in(T_I,T_I+T_0],
   \end{array}\right.\quad\mbox{where}\quad I>0.
\end{equation}
The equilibrium of (\ref{eq1}) with $I_{in}(t)=I$ is given by
\begin{equation}\label{vIhI}
  \left(\begin{array}{c}
    v_I\\ h_I\end{array}\right)=\frac{I}{g_po_h+mg_h}\left(\begin{array}{c} o_h\\ -m\end{array}\right)
\end{equation}
and, for any $k\ge 0,$ the constant $k_i$ of (\ref{ki}) doesn't depend on $i$ and  computes as
\begin{equation}\label{bark}
   \bar k=\left(\sqrt{k}+\sqrt{mv_I^2+g_h h_I^2}\right)^2.
\end{equation}


\begin{cor}\label{propneuron} Let 
$
   g_p,g_h,m,o_h>0
$
and consider $k>0.$ Assume that
\begin{equation}\label{tdneuron}
  \min\{T_I,T_0\}\ge \frac{1}{2\min\{g_p,o_h\}}\ln\frac{\bar k}{k}=:\tau_d.
\end{equation}
Then, the dynamics of any solution $t\mapsto (v(t),h(t))$ of (\ref{eq1}) with the control function (\ref{controlI}) and with the initial condition $(v(0),h(0))=0$ satisfies
\begin{eqnarray}
 &  mv(t)^2+g_hh(t)^2\le k, &  t=(T_I+T_0)\cdot j,\quad j\in\mathbb{N}.\label{dynamics1}\\
   &m(v(t)-v_I)^2+g_h(h(t)-h_I)^2\le k,& t=(T_I+T_0)\cdot j + T_I,\quad j\in\mathbb{N},\label{dynamics2}\\
  & v(t)\le v_I+\sqrt{\dfrac{\bar k}{m}},& t\ge 0.\label{dynamics3}
\end{eqnarray}
In particular, the neuron model (\ref{eq1})-(\ref{reset}) exhibits just sub-threshold oscillations (never develops spiking), if 
\begin{equation}\label{vth}
   v_{th}> v_I+\sqrt{\dfrac{\bar k}{m}}.
\end{equation}
\end{cor}

\noindent {\bf Proof.} The conclusions (\ref{dynamics1})-(\ref{dynamics2}) are direct consequences of (\ref{dyn1}) of Theorem~\ref{proplinear} and we only  need to  explain how (\ref{dyn2}) implies (\ref{dynamics3}). Indeed, (\ref{dyn2}) literally says
\begin{eqnarray*}
  & m(v(t)-v_I)^2+g_h(h(t)-h_I)^2\le \bar k, & t\in\left[(T_I+T_0)j,(T_I+T_0)j+T_I\right],\quad j\in\mathbb{N},\\
   & mv(t)^2+g_hh(t)^2\le \bar k, & t\in\left[(T_I+T_0)j+T_I,(T_I+T_0)(j+1)\right],\quad j\in\mathbb{N},
\end{eqnarray*}
which implies (\ref{dynamics3}) because $v_I>0$ by (\ref{controlI}) and (\ref{vIhI}).\qed

\vskip0.2cm


\begin{figure}[h]\center
 \hskip-0.0cm\includegraphics[scale=0.55]{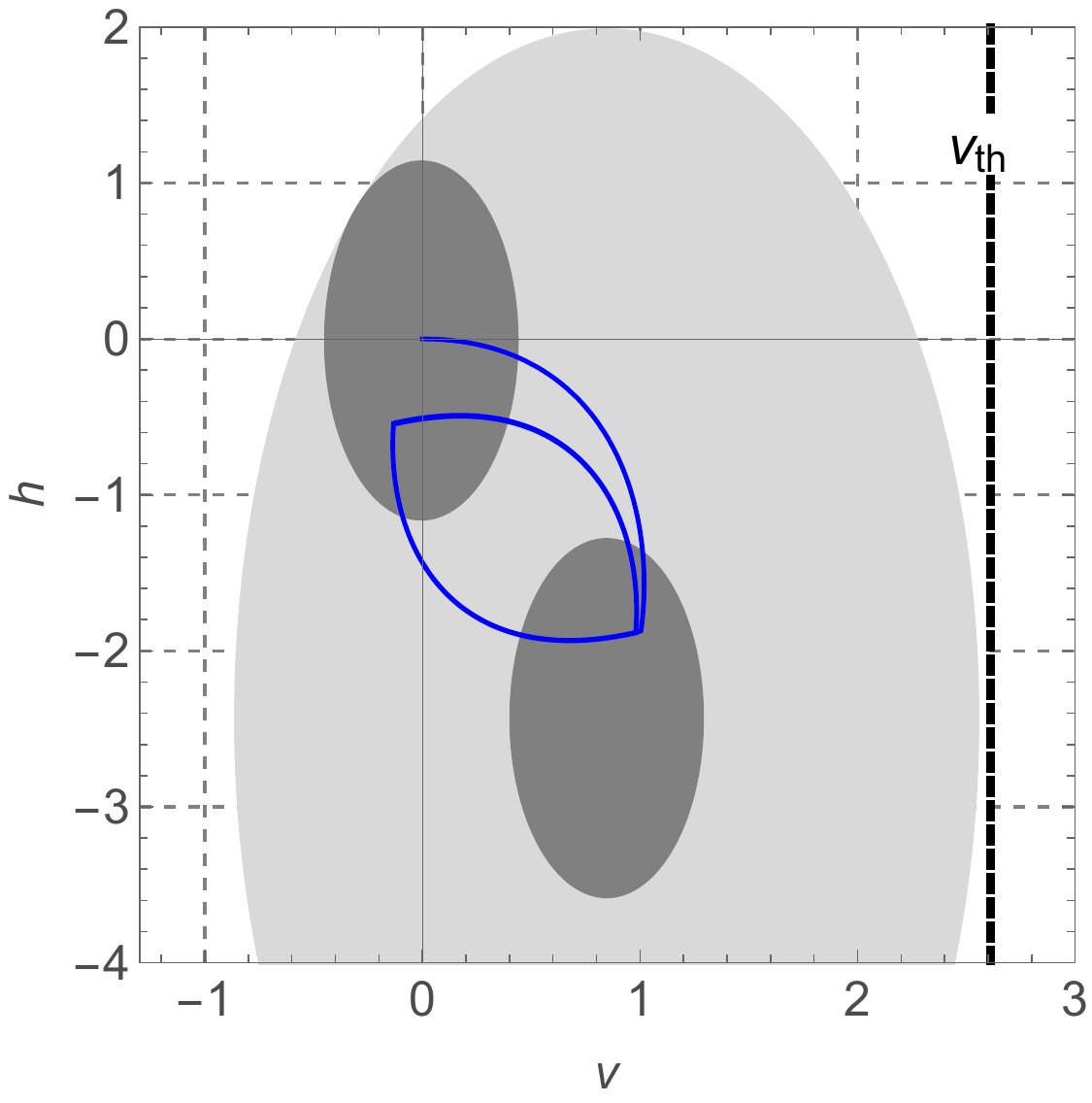}\qquad
\includegraphics[scale=0.55]{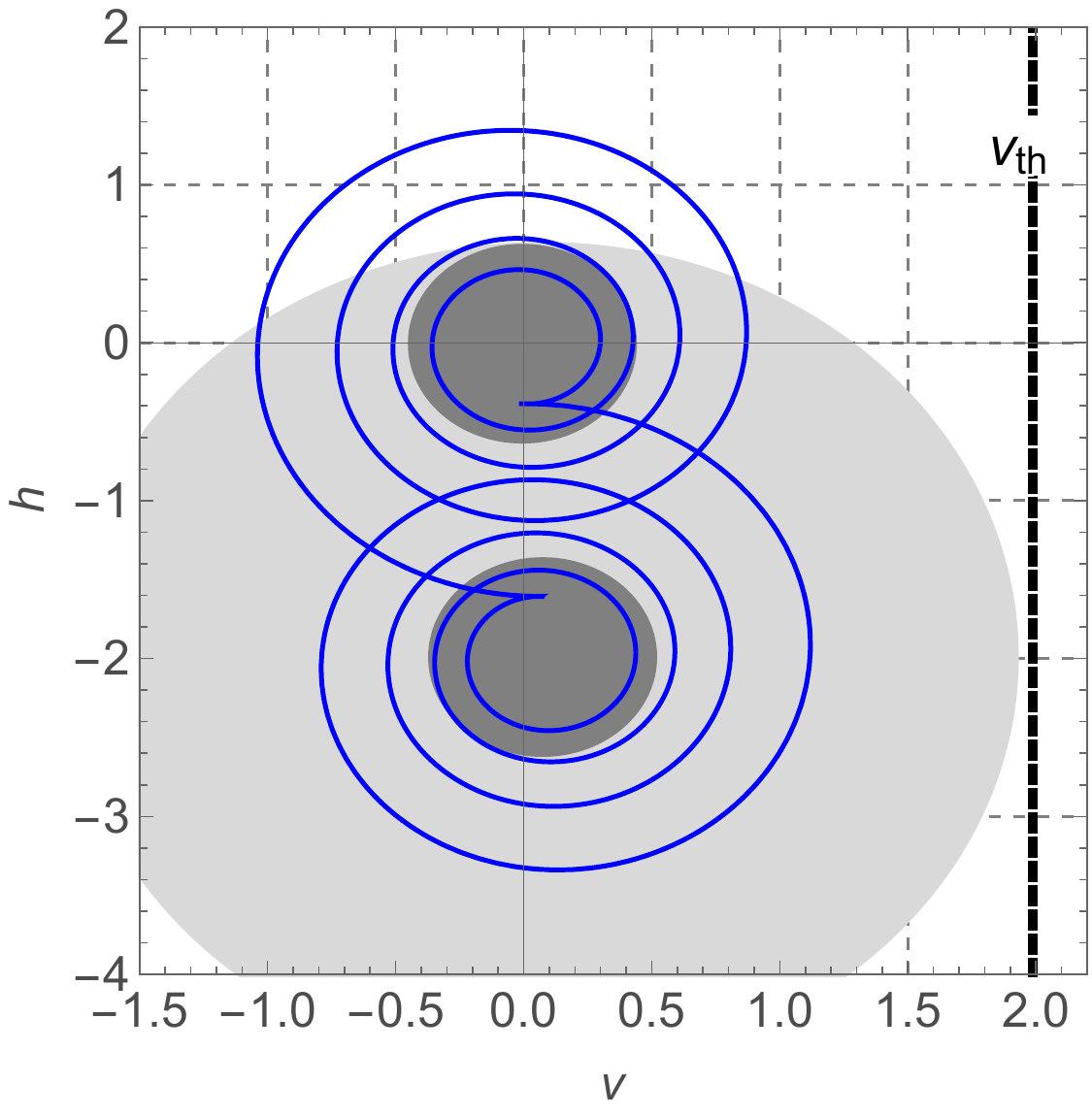}
\vskip-0.4cm
\caption{\footnotesize Left: The trajectory $t\mapsto(v(t),h(t))$ of system (\ref{eq1}) with the initial condition $(v(0),h(0))=0$ for the parameters of  $g_p=0.75$, $g_h=0.15$, $m=1,$ $o_h=0.35$ (from Hasselmo-Shay \cite{hasselmo}) and with the input $I_{in}(t)$ alternating between $0$ and $I=1$ every $T=3.84$ units of time (i.e. $T_0=T_I=3.84$). Right: The attractor of system (\ref{eq1}) for the parameters $g_p=0.04$, $g_h=0.5$, $m=1,$ $o_h=0.04$, whose input $I_{in}(t)$ alternates between $0$ and $I=1$ with period $T_0=T_I=35.7.$ In both figures the dark gray disks are  $N_0^k$ and $N_I^k$, $k=0.2,$ and the light disk is $N_I^{\bar \eps}$, see (\ref{bark}). The line $v=v_th$ is an example of firing threshold that doesn't cause spiking (because the line $v=v_{th}$ does intersect $N_I^{\bar \eps}$).
} \label{figure12}
\end{figure}

\begin{figure}[h]\center
 \hskip-0.0cm\includegraphics[scale=0.55]{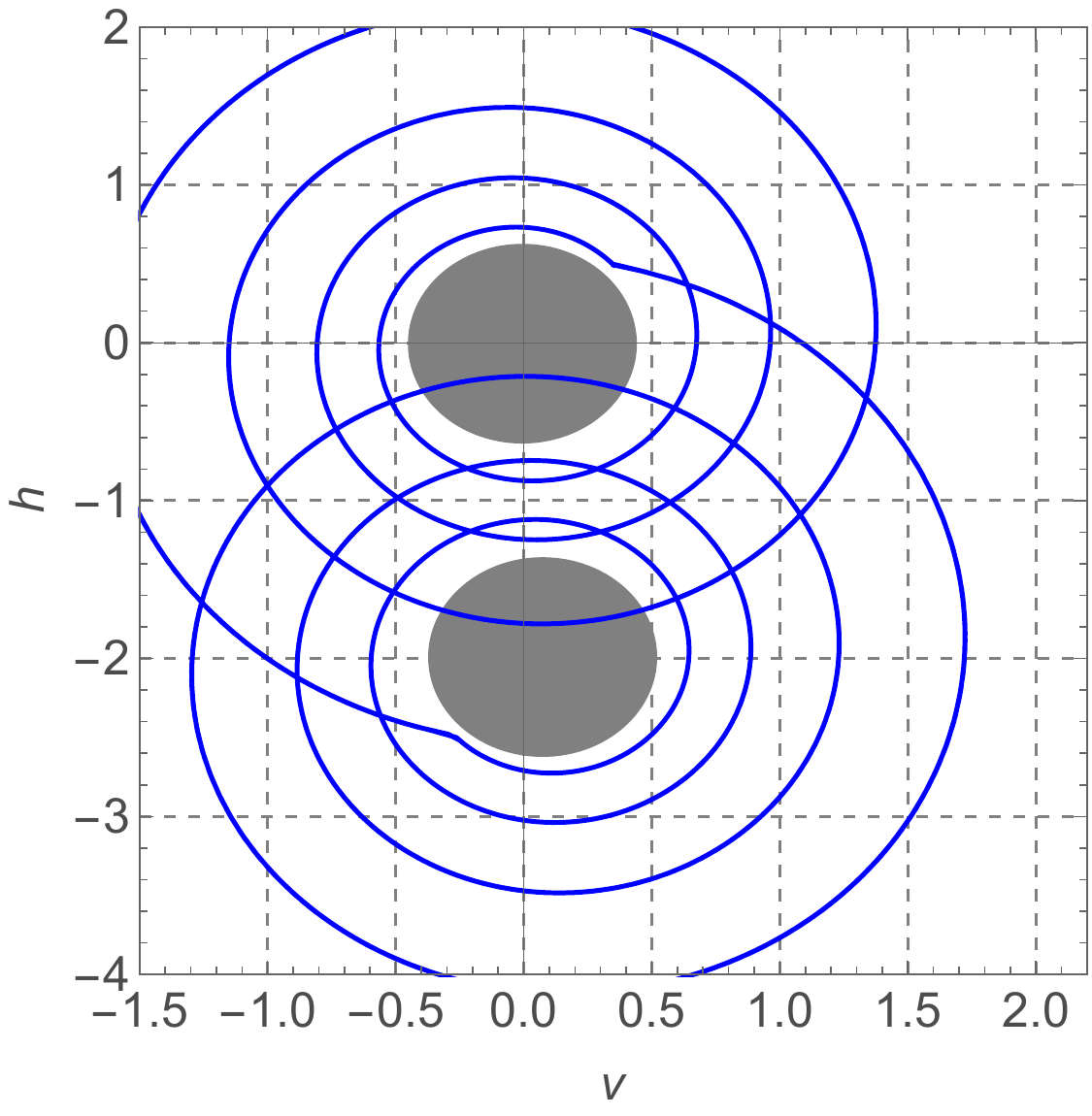}\qquad
\includegraphics[scale=0.55]{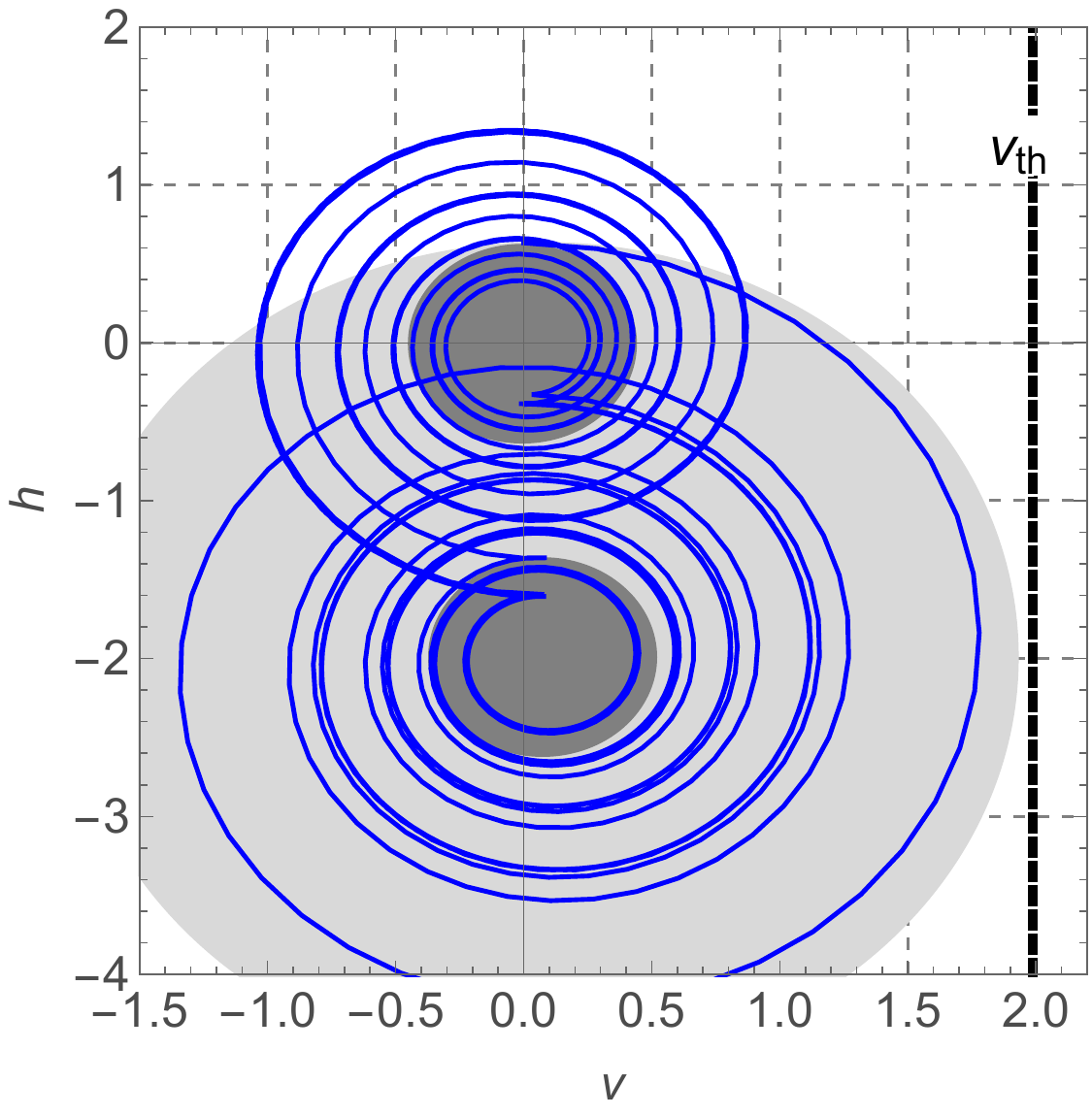}
\vskip-0.4cm
\caption{\footnotesize Both the figures are plotted with the parameters of Fig.~\ref{figure12}(right) except for $T_0$ and $T_I$. Left: Attractor of (\ref{eq1}) for $T_0=T_I=32.$ Right: The solution of (\ref{eq1}) with the initial condition at the top of $N_0^k$ and $T_0=T_I=35.7.$ The meaning of gray disks as well as that of $k$ and $\bar k$ is the same as in Fig.~\ref{figure12}.} \label{figure34}
\end{figure}

\begin{rem} It is possible to develop an analogue of Corollary~\ref{propneuron} (along with the respective analogue of Theorem~\ref{proplinear}), where an estimate on $T_I+T_0$ replaces the estimate (\ref{tdneuron}) on $\min\{T_I,T_0\}.$ The former estimate is known as an {\it average dwell time}, see Liberzon \cite{liberzon}, Yin et al \cite{yin}.
\end{rem}

\noindent Simulations of Figs.~\ref{figure12}-\ref{figure34} illustrate the accuracy of the predictions of Corollary~\ref{propneuron}. At Fig.~\ref{figure12}(left) we drew the solution of the linear neuron model (\ref{eq1}) with the parameters of Hasselmo-Shay \cite{hasselmo} ($g_p=0.75$, $g_h=0.15$, $m=1,$ $o_h=0.35$), $I=1,$ $k=0.2$ and the periods $T_0=T_I=3.84$,  that was computed using the dwell-time formula (\ref{tdneuron}) (which returned the value of $\tau_d=3.836$). Formula (\ref{vth}) provides the estimate $v_{th}>2.56$ for the firing threshold to ensure non-spiking. A possible firing threshold $v_{th}$ is drawn in Fig.~\ref{figure12}(left). The figure also illustrates the construction beyond the estimate (\ref{vth}) whose role is to locate the cylinder $N^{\bar k}_I$ to the left from the line $v=v_{th}.$ Fig.~\ref{figure12}(left) is an example where Corollary~\ref{propneuron} leads to a rather conservative estimate for $v_th$. The figure shows that the value $v_{th}$ can actually be much smaller than  $v_{th}=2.56$ (perhaps around $v_{th}=1.3$) for sub-threshold oscillations to not spike. The sharpness of the estimates of Corollary~\ref{propneuron} is seen e.g. with the parameters $g_p=0.04,$ $g_h=0.5,$ $m=1$, $o_h=0.04$, $k=0.2,$ $I=1$. The dwell time $\tau_d$ given by Corollary~\ref{propneuron} is now $\tau_d=35.621,$ which was used in simulations of Fig.~\ref{figure12}(right) (we took $T_0=T_I=35.7$) where  the respective attractor of model (\ref{eq1}) is shown. First of all, one can see that the switchings (corners of the trajectory) occur very close to the boundary of the cylinders $N_0^k$ and $N_I^k$. Moreover, Fig.~\ref{figure34}(left) shows that the switching points are no longer in $N_0^k$ and $N_I^k,$ if $T_0$ and $T_I$ reduce to $T_0=T_I=32$, which confirms that $\tau_d=35.621$ is a relatively sharp dwell time bound. Finally, Fig.~\ref{figure34}(d) illustrates that the estimate (\ref{vth}) for the maximal current is also accurate, i.e. a trajectory with the initial condition in $N_0^k$ can pass quite close to the rightmost point of the ellipse $N_I^{\bar k}.$

\section{An extension in the multi-dimensional nonlinear case}

\begin{figure}[h]\center
\includegraphics[scale=0.55]{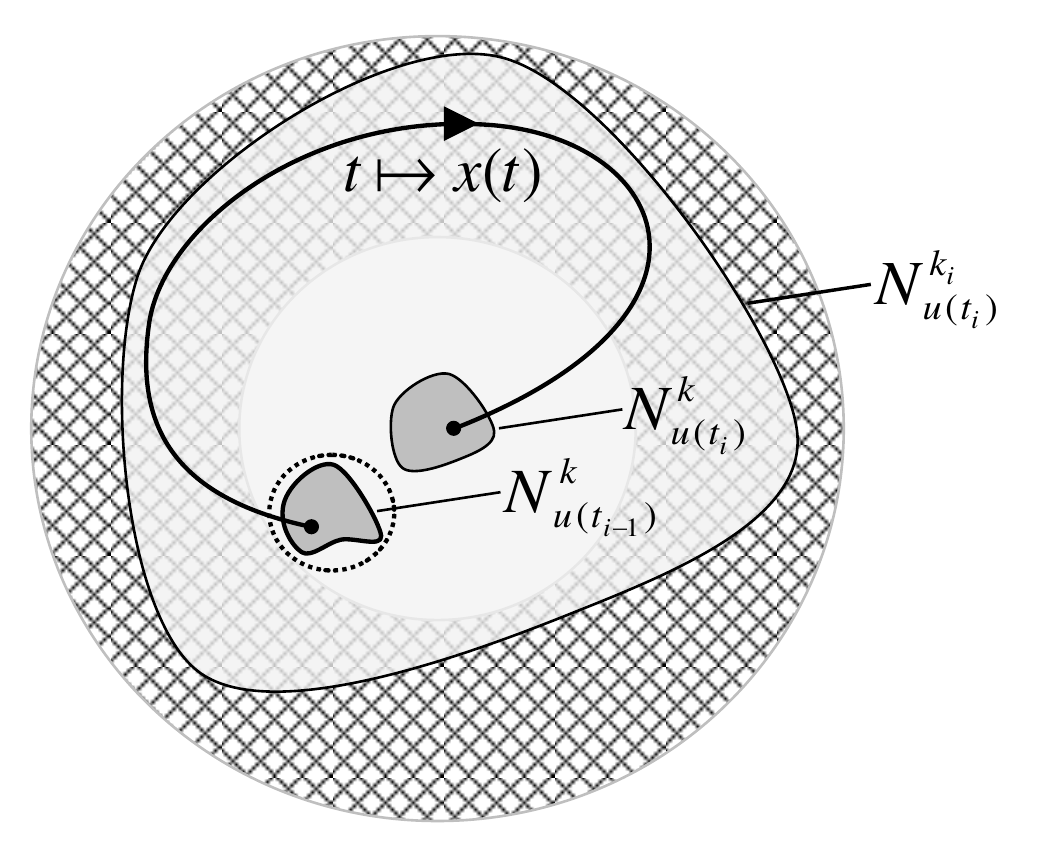}
\vskip-0.2cm
\caption{\footnotesize  The relative location of the level curves  $\partial N_{u(t_{i})}^{k_i}$ and $\partial N_{u(t_{i-1})}^{k}$ along with a solution of switched system (\ref{ss}) with initial condition $x(t_{i-1})\in N_{u(t_{i-1})}^k$ on the interval $[t_{i-1},t_i].$ The dotted circle is the circle of radius $\alpha^{-1}_{u(t_{i-1})}$ centered at $x_{u(t_{i-1})}$ and which, by (\ref{ab}), surrounds $N_{u(t_{i-1})}^{k}$. The textured ring is the ring centered at $x_{u(t_i)}$ with the inner radius $\|x_{u(t_i)}-x_{u(t_{i-1})}\|+\alpha^{-1}_{u(t_{i-1})}(k)$ (i.e. the minimal radius for which the ring surrounds the dotted circle) and with the outer radius $\alpha_{u(t_i)}^{-1}(k_i)$, so that $\partial N_{u(t_{i})}^{k_i}$ is contained in this ring by (\ref{ab}).} \label{figXX}
\end{figure}

\noindent When (\ref{ss}) is nonlinear and $V_u$ is an arbitrary Lyapunov function, closed-form formulas for the dwell-time can be obtained when $V_u$ admits the estimates  
\begin{eqnarray}
&&   \alpha_u(\|x-x_u\|)\le V_u(x)\le \beta_u(\|x-x_u\|),\qquad x\in\mathbb{R}^n,\label{ab}\\
&& (V_u)'(x)f_u(x)\le -\eps_u V_u(x),\qquad x\in\mathbb{R}^n,\label{k}
\end{eqnarray}
where $\alpha$, $\beta$ are strictly monotonically increasing functions with $\alpha_u(0)=\beta_u(0)$, and $\eps_u>0.$ Conditions (\ref{ab}) does appear in Alpcan-Basar \cite{alpcan}, but it is not explode in \cite{alpcan} for computing the dwell time.
Based upon Makarenkov-Phung \cite{anthony}, we can offer the following computational formulas to estimate the location of the the dynamics of switched system (\ref{ss}) (see Fig.~\ref{figXX} for explanation of the crucial constant (\ref{ki1})).

\begin{thm}\label{nonlinear} Assume that, for any $u\in\mathbb{R}$,  system (\ref{ssi})
 with $f_u\in C^1(\mathbb{R}^n,\mathbb{R}^n)$,
 admits an equilibrium $x_u$ whose Lyapunov function $V_u\in C^1(\mathbb{R}^n,\mathbb{R})$ satisfies (\ref{ab})-(\ref{k}) with strictly increasing  $\alpha_u,\beta_u\in C^0(\mathbb{R},\mathbb{R})$ satisfying $\alpha_u(0)=\beta_u(0)=0$ and $\eps_u>0.$
If the successive discontinuities $t_1,t_2,...$ of the control signal $u(t)$ verify
 \begin{equation}\label{sat2new1}
   t_i-t_{i-1}\ge -\frac{1}{\eps_{u(t_i)}}\ln\frac{k}{k_i},
 \end{equation}
 where 
 \begin{equation}\label{ki1}
 k_i=\beta_{u(t_i)}\left(\|x_{u(t_i)}-x_{u(t_{i-1})}\|+\alpha_{u(t_{i-1})}^{-1}(k)\right),\quad i\in\mathbb{N},
 \end{equation}
 then
$$
   x(t_i)\in N_{u(t_i)}^k,\quad i\in\mathbb{N},
$$
and
\begin{equation}\label{propro}
   x(t)\in N^{k_i}_{u(t_i)},\quad t\in[t_{i-1},t_i],\ i\in\mathbb{N}.
\end{equation}
\end{thm}

\vskip0.2cm

\noindent The crucial difference between Theorems~\ref{proplinear} and \ref{nonlinear} is seen from Figs.~\ref{figX} and \ref{figXX}. Indeed, the level set $N_{u(t_i)}^{k_i}$ used in Theorem~\ref{proplinear} is the minimal possible level set that contains $N^k_{u(t_{i-1})}$ ($\partial N_{u(t_i)}^{k_i}$ just touches $\partial N^k_{u(t_{i-1})}$ in Fig.~\ref{figX}). In contrast, the way how we define $N_{u(t_i)}^{k_i}$ in Theorem~\ref{proplinear} ($\partial N_{u(t_i)}^{k_i}$ is inscribed into a ring surrounding $\partial N^k_{u(t_{i-1})}$) is more conservative as it may leave a significant gap between $\partial N^k_{u(t_{i-1})}$ and $N_{u(t_i)}^{k_i}$ (that is seen in Fig.~\ref{figXX}).

\vskip0.2cm

\noindent Theorem~\ref{nonlinear} appears in our other manuscript \cite{anthony}, but without (\ref{propro}) which is the crucial quantity when comparing Theorems~\ref{nonlinear} and \ref{proplinear} (see Fig.~\ref{figXX}).

\vskip0.2cm

\noindent {\bf Proof.} Let us fix $i\in \mathbb{N}$ and consider a solution $x$ of (\ref{ss}) with $x(t_{i-1})\in N_{u_{i-1}}^k.$ Our goal is to show that $x(t_i)\in N_{u_i}^k.$
Given $k>0$, define $k_i>0$ according to (\ref{ki}). 
  As in the proof of Theorem~\ref{proplinear}, introduce 
$
   v(t)=V_{u_i}(x(t)).
$
 By construction, see illustration in  Fig.~\ref{figXX}, 
$N_{u_i}^{k_i}\supset N_{u_{i-1}}^k$ and so 
$x(t_{i-1})\in N_{u_i}^{k_i}$.  On the other hand,  by (\ref{k}) we have
$$  
\dot v(t)\le -\eps_{u(t_i)} v(t),\quad t_{i-1}\le t\le t_i, \quad i\in\mathbb{N}.
$$
 Therefore, using (\ref{sat2new1}), we obtain
$$
  v(t_i)=e^{-\eps_{u(t_i)}(t_i-t_{i-1})}k_i \le k,
$$
which completes the proof.\qed

\vskip0.2cm


\section{Conclusion} \noindent In this paper we offered sharp dwell time formulas for the frequency of the successive switchings of a planar switched affine system, that traps the solution in a given tube that connects the successive equilibria of individual subsystems. This is the first paper where the respective estimate for the location of the dynamics of switched systems is used in the context of neuroscience. Specifically, we gave explicit condition for a linear neuron model to never reach the firing threshold, i.e. to operate in just subthreshold mode. Since non-spiking in neuron model represents the main motivation for this paper, we focused on the situation where the homogeneous part of the linear subsystem stays constant and the switching occurs in the inhomogeneous part only. However, we explained (Remark~\ref{rem1}) how the results of the paper extends to switching between arbitrary affine systems. Furthermore, we considered 2-dimensional affine systems only, but we don't see any obstacles for the extension of the analysis to the multi-dimensional case. It also looks doable to account for possible uncertainties in switched system (\ref{ss}) complementing the global results of Jin et al \cite{jin}.  
 The ideas of the paper can be further used in power electronics where planar switched affine systems model  power converters \cite{converters1,converters2,converters3}. 
 
 \vskip0.2cm
 
\noindent On top of the above, we built upon \cite{anthony} and offered (Section~3) explicit formulas to estimate the trapping region of the dynamics of (\ref{ss}) in the case or arbitrary Lyapunov functions (not necessary quadratic), which might show its effectiveness in combination with polynomial Lyapunov functions of \cite{po2,po3,po4,po5}.

\vskip0.2cm

\noindent Our most immediate upcoming plans include developing the ideas of the paper to the level capable to make contributions in local stability of switched genetic regulatory networks \cite{zhang,yao}) and multi-agent systems \cite{gao,shucker}.

\section{Acknowledgments.} \noindent The first author is  supported by NSF Grant CMMI-1436856 and by Burroughs Wellcome Fund Collaborative Research Travel Grant \#1017453.

\bibliographystyle{plain}

\section{References}
\begin{thebibliography}{000}


\bibitem{alpcan} T. Alpcan, T. Basar, 
A stability result for switched systems with multiple equilibria, {
Dyn. Contin. Discrete Impuls. Syst. Ser. A Math. Anal.} {\bf 17}, (2010), no. 6, 949--958. 


\bibitem{po2} M. A. Ben Sassi, S. Sankaranarayanan, X. Chen, E. Abraham, 
Linear relaxations of polynomial positivity for polynomial Lyapunov function synthesis.  
IMA J. Math. Control Inform. 33 (2016), no. 3, 723--756. 

\bibitem{bolzen}
P. Bolzern, W. Spinelli, Quadratic stabilization of a switched affine system about a nonequilibrium point, Proceeding of the 2004 American Control Conference, June 30--July 2, 2004, 3890--3895.

\bibitem{non-spiking2} W. Chen, R. Maex, R. Adams, V. Steuber, L. Calcraft, N. Davey, Clustering predicts memory performance in networks of
spiking and non-spiking neurons, Frontiers in Computational Neuroscience 5 (2011), no. 5, Article 14, 10pp.

\bibitem{coombes} S. Coombes, R. Thul, K. C. A. Wedgwood, 
Nonsmooth dynamics in spiking neuron models.  
Phys. D 241 (2012), no. 22, 2042--2057. 

\bibitem{dorothy}
M. Dorothy, S.-J. Chung, 
Switched systems with multiple invariant sets.  
Systems Control Lett. 96 (2016), 103--109. 


\bibitem{converters1}
P. Gupta and A. Patra, Hybrid mode switched control of DC--DC boost
converter circuits, IEEE Trans. Circuits Syst. II, Exp. Briefs 52 (2005), no.  11, 734--738.

\bibitem{jin} Y. Jin, J. Fu, Y. Zhang, Y. Jing, 
Reliable stabilization of switched system with average dwell-time approach. 
J. Franklin Inst. 350 (2013), no. 3, 452--463. 

\bibitem{gao}
L. Gao, Y. Cui, X. Xu, Y. Zhao, 
Distributed consensus protocol for leader-following multi-agent systems with functional observers.  
J. Franklin Inst. 352 (2015), no. 11, 5173--5190. 

\bibitem{hasselmo}
M. E. Hasselmo, C. F. Shay, Grid cell firing patterns may arise from feedback interaction
between intrinsic rebound spiking and transverse traveling
waves with multiple heading angles, Frontiers in Systems Neuroscience 8 (2014), no. 8, Article 201, 24~pp.

\bibitem{izh}
Izhikevich, E. M. (2007) Dynamical systems in neuroscience: the geometry of excitability and bursting. Computational Neuroscience. MIT Press, Cambridge, MA, xvi+441 pp. 

\bibitem{78} H. K. Khalil, Nonlinear systems. Macmillan Publishing Company, New York, 1992.
xii+564 pp.

\bibitem{li} J. Li, R. Ma, G. M. Dimirovski, J. Fu, 
Dwell-time-based stabilization of switched linear singular systems with all unstable-mode subsystems. 
J. Franklin Inst. 354 (2017), no. 7, 2712--2724. 

\bibitem{liberzon} D.
Liberzon, {\rm Switching in systems and control.} Systems \& Control: Foundations \& Applications. Birkhauser Boston, Inc., Boston, MA, 2003. 


\bibitem{anthony} O. Makarenkov, A. Phung  Dwell time for switched systems with multiple equilibria on a finite time-interval, submitted. arXiv:1703.06205.

\bibitem{po3}
L. Menini, A. Tornambe, 
On a Lyapunov equation for polynomial continuous-time systems.
Internat. J. Control 87 (2014), no. 2, 393--403. 

\bibitem{po4}
M. Peet,
Exponentially stable nonlinear systems have polynomial Lyapunov functions on bounded regions. (English summary) 
IEEE Trans. Automat. Control 54 (2009), no. 5, 979--987. 


\bibitem{fridman} A. Polyakov, L. Fridman, Stability notions and Lyapunov functions for sliding mode control systems. J. Franklin Inst. 351 (2014), no. 4, 1831--1865.


\bibitem{converters2}
 A. Schild, J. Lunze, J. Krupar, and W. Schwarz,
Design of generalized hysteresis controllers for dc-dc
switching power converters, IEEE Trans. Power
Electron. 24 (2009), no. 1, 138--146.

\bibitem{shucker}
B. Shucker, T. D. Murphey, J. K. Bennett, Convergence-Preserving Switching for Topology-Dependent Decentralized Systems, IEEE Transactions on Robotics 24 (2008) 1405--1415.

\bibitem{converters3}
P. Siewniak, B. Grzesik, The piecewise-affine model of buck converter suitable for practical
stability analysis, Int. J. Circ. Theor. Appl. 43 (2015) 3--21.

\bibitem{po5}
W. Tian, C. Zhang, C. Qian, S. Li, Global stabilization of inherently non-linear systems using continuously differentiable controllers.
Nonlinear Dynam. 77 (2014), no. 3, 739--752. 

\bibitem{non-spiking1}
C. Vich, A. Guillamon, 
Dissecting estimation of conductances in subthreshold regimes. 
J. Comput. Neurosci. 39 (2015), no. 3, 271--287. 

\bibitem{77} M. Vidyasagar, Nonlinear systems analysis. Reprint of the second (1993) edition.
Classics in Applied Mathematics, 42. Society for Industrial and Applied Mathematics
(SIAM), Philadelphia, PA, 2002. xviii+498 pp.

\bibitem{xu}
H. Xu, Y. Zhang, J. Yang, G. Zhou, L. Caccetta, Practical exponential set stabilization for switched nonlinear systems with multiple subsystem equilibria, { J. Global Optim.} {\bf 65} (2016), no. 1, 109--118.

\bibitem{yao} Y. Yao, J. Liang, J. Cao, 
Stability analysis for switched genetic regulatory networks: an average dwell time approach. 
J. Franklin Inst. 348 (2011), no. 10, 2718--2733. 


\bibitem{yin}
Y. Yin, G. Zong, X. Zhao, 
Improved stability criteria for switched positive linear systems with average dwell time switching. 
J. Franklin Inst. 354 (2017), no. 8, 3472--3484. 

\bibitem{michel}
G. Zhai, B. Hu, K. Yasuda, A. Michel, Stability analysis of switched systems with
stable and unstable subsystems: an average dwell time approach. Internat. J. Systems
Sci. 32 (2001), no. 8, 1055--1061.

\bibitem{zhang}
W. Zhang, J. Fang, W. Cui, 
Exponential stability of switched genetic regulatory networks with both stable and unstable subsystems. 
J. Franklin Inst. 350 (2013), no. 8, 2322--2333. 

\bibitem{automatica} Y. Zhang, O. Makarenkov, N. Gans, Extremum seeking control of a nonholonomic system with sensor constraints. Automatica J. IFAC 70 (2016) 86--93.



\end{thebibliography}

\end{document}